\documentclass{article}
\usepackage{amsmath,amsfonts,amssymb,amsthm}
\title{Modified scattering for a wave equation with weak dissipation}
\author{Jens Wirth\thanks{Department of Applied Analysis, TU Bergakademie Freiberg, 09596 Freiberg, Germany}\ \thanks{email: {\tt wirth@math.tu-freiberg.de}}}
\def\R{\mathbb R}
\def\d{\mathrm d}\def\D{\mathrm D}
\def\H{\mathcal H}\def\J{\mathcal J}
\DeclareMathOperator*{\slim}{s-lim}
\DeclareMathOperator{\supp}{supp}
\newtheorem{thm}{Theorem}[section]
\newtheorem{lem}[thm]{Lemma}
\newtheorem{prop}[thm]{Proposition}
\newtheorem{cor}[thm]{Corollary}
\numberwithin{equation}{section}
\begin{document}
\maketitle
\begin{abstract}
We consider the Cauchy problem for the weakly dissipative wave equation
$$ \square u+\frac\mu{1+t} u_t=0 $$
with parameter $\mu\ge2$.

Based on the explicit representations of solutions provided in [Math. Meth.
Appl. Sci. 2004; {\bf 27}:101-124] sharp decay estimates for data from a 
dense subspace of the energy space are derived. Furthermore, sharpness
is discussed in terms of a modified scattering theory.

\end{abstract}
\section{Introduction}
We are interested in a precise description of the behaviour of solutions to the Cauchy problem
\begin{equation}\label{eq:CP}
  \square u+\frac{\mu}{1+t} u_t=0,
  \qquad\qquad u(0,\cdot)=u_1,\quad u_t(0,\cdot)=u_2
\end{equation}
with $\square=\partial_t^2-\Delta$ and for data $u_1\in H^1(\R^n)$, $u_2\in L^2(\R^n)$. It is well known that this problem is well-posed in the sense that it has a unique solution in the space 
$$ C^1([0,\infty),L^2(\R^n))\cap C([0,\infty),H^1(\R_+)). $$

It is known from the fundamental works of A.Matsumura, \cite{Mat77}, and
H.Uesaka, \cite{Ues79}, that the energy of the solutions to this equation,
\begin{equation}
  E(u;t)=\frac12\int (u_t^2+|\nabla u|^2)\d x,
\end{equation}
satisfies an estimate of the form
\begin{equation}
  E(u;t)=\mathcal O(t^{-\alpha}),\qquad \alpha=\min\{2,\mu\}.
\end{equation}
Following F.Hirosawa and H.Nakazawa, \cite{HN04}, it is possible to improve the last estimate in the case $\mu>2$ to
\begin{equation}
  E(u;t)=o(t^{-2}),
\end{equation}
which means that the above given estimate is not sharp (with respect to the strong
topology).

In \cite{Wir04} the author provided explicit representations of solutions to the 
Cauchy problem \eqref{eq:CP} in terms of special functions. Furthermore, from that
treatment it follows that the estimate $E(u;t)=\mathcal O(t^{-2})$ is sharp in the sense of a norm estimate for the energy operator 
\begin{equation}
 \mathbb E(t) : (\langle \D\rangle u_1,u_2)^T\mapsto (|\D|u(t,\cdot),\partial_t u(t,\cdot))^T 
\end{equation}
associated to the solution representation. As usual we denote $\langle\xi\rangle 
=\sqrt{1+|\xi|^2}$.

In this paper we will use this explicit representation to describe a dense subspace of the energy space such that for all data from that space we obtain a sharp two-sided estimate of the energy in the form
\begin{equation}
  E(u;t)\sim t^{-\mu}. 
\end{equation}
The result can be formulated in terms of a modified scattering theory. Especially
it describes the discrepancy between the estimates of Hirosawa/Nakazawa and the
norm estimate from \cite{Wir04}.

\section{Preliminaries}
Let us recall the representations of solutions to equation \eqref{eq:CP}. With the notation
\begin{subequations}
\begin{align}
  \Psi_{k,s,\rho,\delta}(t,\xi) &= |\xi|^k\langle\xi\rangle^{s+1-k}
  \begin{vmatrix}\H^-_\rho(|\xi|)&\H_{\rho+\delta}^-((1+t)|\xi|)\\
  \H_\rho^+(|\xi|)&\H_{\rho+\delta}^+((1+t)|\xi)\end{vmatrix}\\
  &=2i\csc(\rho\pi)|\xi|^k\langle\xi\rangle^{s+1-k}\begin{vmatrix}
  \J_{-\rho}(|\xi|)&\J_{-\rho-\delta}((1+t)|\xi)\\
  (-)^\delta \J_\rho(|\xi|)& \J_{\rho+\delta}((1+t)|\xi|)\end{vmatrix},
\end{align}
\end{subequations}
the last line analytically continued to $\rho\in\mathbb Z$, the following theorem is valid, \cite[Theorem 2.1]{Wir04}.

\begin{thm}
  The spatial Fourier transform of the solution to problem \eqref{eq:CP} can 
  be represented in the form
  \begin{equation} 
        \hat u(t,\xi)=\sum_{j=1,2} \Phi_j(t,\xi)\hat u_j(\xi) ,
  \end{equation}
  where
  \begin{subequations}
  \begin{align}
    \Phi_1(t,\xi)&=\frac{i\pi}4 (1+t)^\rho \Psi_{1,0,\rho-1,1}(t,\xi),\\
    \partial_t \Phi_1(t,\xi)&=\frac{i\pi}4 \Psi_{2,1,\rho-1,0}(t,\xi),\\
    \Phi_2(t,\xi)&=-\frac{i\pi}4 (1+t)^\rho \Psi_{0,-1,\rho,0}(t,\xi),\\
    \partial_t \Phi_2(t,\xi)&=-\frac{i\pi}4 (1+t)^\rho \Psi_{1,0,\rho-1,0}(t,\xi)
  \end{align}
  \end{subequations}
  with $\rho=(1-\mu)/2$.
\end{thm}

We want to obtain estimates in $L^2$-scale, which correspond by Plancherel's theorem to $L^\infty$ estimates for the Fourier multiplier. Therefore we recall the fundamental estimate from \cite{Wir04}.

\begin{lem}\label{lem2.2}
 It holds $\Psi_{k,s,\rho,\delta}(t,\cdot)\in L^\infty(\R^n)$ for all $t$ if and only if $s\le0$ and $k\ge|\delta|$. Furthermore, the estimate
\begin{equation}
  ||\Psi_{k,s,\rho,\delta}(t,\cdot)||_\infty\sim 
  \begin{cases} (1+t)^{-\frac12},& \rho\ne0, |\rho|-k\le-\frac12,\\
  (1+t)^{|\rho|-k},&\rho\ne0, |\rho|-k\ge-\frac12,\\
  (1+t)^{-k} \log(e+t),&\rho=0, k\le\frac12 \end{cases}
\end{equation}
is valid.
\end{lem}

For later use we will introduce a notation. Let 
\begin{equation}
  [\xi]=\frac{|\xi|}{\langle\xi\rangle}.
\end{equation}
By the aid of this symbol we can control the vanishing order in the frequency
$\xi=0$. We denote
\begin{equation}
 [\D]^\kappa L^2(\R^n)=\{ [\D]^\kappa f\;|\; f\in L^2\},\qquad ||g||_{[\D]^\kappa L^2}=||[\D]^{-\kappa}g||_2. \end{equation}
It turns out that the usage of date from these spaces allows an improvement of the
decay order of the energy.

\section{Improvements of the energy decay}

\begin{thm}\label{thm:3.1} 
Let $\langle\D\rangle u_1,u_2\in [\D]^\kappa L^2(\R^n)$. Then the solution
$u=u(t,x)$ to \eqref{eq:CP} with $\mu>2$ satisfies 
\begin{equation}
  || |\D|u(t,\cdot)||_2+||u_t(t,\cdot)||_2\lesssim 
    || \langle\D\rangle u_1,u_2||_{[\D]^\kappa L^2}
    \begin{cases} (1+t)^{-1},&\kappa=0,\\
     (1+t)^{-1-\kappa},&0\le\kappa\le\frac{\mu-2}2,\\
     (1+t)^{-\frac\mu2},&\kappa\ge\frac{\mu-2}2. \end{cases}
\end{equation}
\end{thm}
\begin{proof}
The proof is a direct application of Lemma~\ref{lem2.2}. We use the canonical 
isomorphism between $[\D]^\kappa L^2$ and $L^2$ to simplify the estimates.
The mapping $(\langle\D\rangle u_1,u_2)^T\mapsto (|\D|u(t,\cdot),u_t(t,\cdot)^T$ can be represented as matrix Fourier multiplier in terms of $\Psi_{k,s,\rho,\delta}$,
because it includes the identification of spaces we include the number $\kappa$ in the notation. It holds
$$ \mathbb E^{(\kappa)}(t,\xi)=\frac{i\pi}4(1+t)^{\rho}
\begin{pmatrix}
  \Psi_{2+\kappa,0,\rho-1,1}(t,\xi) & -\Psi_{1+\kappa,0,\rho,0}(t,\xi)\\
  \Psi_{2+\kappa,0,\rho-1,0}(t,\xi) & -\Psi_{1+\kappa,0,\rho,-1}(t,\xi)
\end{pmatrix}. $$
It remains to check the conditions on the indices. For $\mu>2$ we have $\rho<-1/2$
and therefore we have to compare $1+\kappa-\rho$ with $-1/2$.  
\end{proof}

Some remarks for the interpretation of the result are necessary. The decay rate
$t^{-\mu/2}$ corresponds to the rate for high frequencies. The cut-off  by $[\xi]^\kappa$ dampes out the small frequencies around the exceptional frequency $\xi=0$. 

For $\mu>2$ the decay rate obtained in \cite{Wir04} was determined by
the small frequencies around $\xi=0$. This is the reason for the improvement of this theorem. An improvement over the decay rate for high frequencies is not possible.

\section{Sharpness}
The aim of this section is to discuss the sharpness of the estimate in the
limit case $\kappa=(\mu-2)/2$. For this we compare the operator family
$\mathbb E^{(\kappa)}(t)$ with the unitary evolution $\mathcal E_0(t)$
of free waves in energy space,
\begin{equation}
 \mathcal E_0(t,\xi)=\begin{pmatrix} \cos t|\xi| & \sin t|\xi|\\
 -\sin t|\xi| & \cos t|\xi| \end{pmatrix},
\end{equation}
modified by the decay rate. 

\begin{thm}\label{thm:4.1}
 Let $\kappa=\frac{\mu-2}2$ and $\mu\ge2$.
  Then the strong limit
  \begin{equation}
   Z_+=\slim_{t\to\infty} (1+t)^{\mu/2} \mathcal E_0(-t)\mathbb E^{(\kappa)}(t) 
   \end{equation}
  exists in the operator space $[\D]^\kappa L^2\to L^2$.
\end{thm}

Let us first explain the main strategy of the proof. Theorem~\ref{thm:3.1} yields a uniform bound for the operator family.
In order to obtain the strong convergence we employ Banach-Steinhaus theorem for the dense subspace
$$ M=\{ f\in L^2\;|\; 0\not\in\supp\hat f\}. $$
For this it is sufficient to prove (uniform) convergence of the Fourier multiplier
for $|\xi|\ge c$ with $c>0$. This can be done by the given explicit representations using known asymptotic expansions for Bessel functions, \cite{WA}.
Because of its simplicity we will use the representation by real valued functions
valid for $\rho\not\in\mathbb Z$. The multiplier $Z_+$ will be analytic in $\rho$
and therefore the statement will follow in all cases.

The result is based on the following asymptotic formula, \cite{WA}.
\begin{prop}\label{pB}
For Bessel functions of first kind it holds
\begin{equation}
  \J_\rho(z)=\sqrt{\frac2{\pi z}}\cos\left(z-\frac{\rho\pi}2-\frac\pi4\right)
  +\mathcal O(z^{-3/2}),\qquad z\to\infty.
\end{equation}
\end{prop}

\begin{proof}[Proof of Theorem~\ref{thm:4.1}] Let $|\xi|\ge c>0$. Then it holds
\begin{align*}
  (1+t)^{\mu/2}\mathbb E^{(\kappa)}(t,\xi)=\frac{i\pi}4 (1+t)^{1/2}
  \begin{pmatrix}\Psi_{2+\kappa,0,\rho-1,1}(t,\xi) & -\Psi_{1+\kappa,0,\rho,0}(t,\xi)\\
  \Psi_{2+\kappa,0,\rho-1,0}(t,\xi) & -\Psi_{1+\kappa,0,\rho,-1}(t,\xi)\end{pmatrix}
\end{align*} 
such that
\begin{align*}
  (1+t)^{\mu/2}& \mathcal E_0(-t,\xi)\mathbb E^{(\kappa)}(t,\xi) \\
  &=\frac{i\pi}4 (1+t)^{1/2}\begin{pmatrix} \cos t|\xi| & -\sin t|\xi|\\
  \sin t|\xi| & \cos t|\xi| \end{pmatrix}\\&\qquad\qquad
  \begin{pmatrix}\Psi_{2+\kappa,0,\rho-1,1}(t,\xi) & -\Psi_{1+\kappa,0,\rho,0}(t,\xi)\\
  \Psi_{2+\kappa,0,\rho-1,0}(t,\xi) & -\Psi_{1+\kappa,0,\rho,-1}(t,\xi)\end{pmatrix}\\
  &= \begin{pmatrix} m_{1,1}(t,\xi)&m_{1,2}(t,\xi)\\m_{2,1}(t,\xi)&m_{2,2}(t,\xi)\end{pmatrix}.
\end{align*}
We proceed with calculating $m_{1,1}(t,\xi)$. It holds
\begin{align*}
 m_{1,1}(t,\xi)=&\frac\pi2\csc(\rho\pi) |\xi|^{2+\kappa}\langle\xi\rangle^{-\kappa-1} (1+t)^{1/2}\\
     &\bigg(\cos t|\xi|\big(\J_{1-\rho}(|\xi|)\J_\rho((1+t)|\xi|)+\J_{\rho-1}(|\xi|)\J_{-\rho}((1+t)|\xi|)\big)\\
     &-\sin t|\xi| \big(\J_{1-\rho}(|\xi|)\J_{\rho-1}((1+t)|\xi|)-\J_{\rho-1}(|\xi|)\J_{1-\rho}((1+t)|\xi|)\big)\bigg)\\
     =&\sqrt{\frac\pi2}\csc(\rho\pi) |\xi|^{2+\kappa-\frac12}\langle\xi\rangle^{-\kappa-1} (1+t)^{1/2}\\
     &\bigg(\cos t|\xi|\big(\J_{1-\rho}(|\xi|)\cos((1+t)|\xi|-\frac{\rho\pi}2-\frac\pi4)\\
     &\qquad\qquad+\J_{\rho-1}(|\xi|)\cos((1+t)|\xi|+\frac{\rho\pi}2-\frac\pi4)\big)\\
     &-\sin t|\xi| \big(\J_{1-\rho}(|\xi|)\cos((1+t)|\xi|-\frac{(\rho-1)\pi}2-\frac\pi4)\\&\qquad\qquad
     -\J_{\rho-1}(|\xi|)\cos((1+t)|\xi|-\frac{(1-\rho)\pi}2-\frac\pi4)\big)\bigg)\\
     &+\mathcal O((1+t)^{-1}|\xi|^{-1})\\
     =& \sqrt{\frac\pi2}\csc(\rho\pi) [\xi]^{\kappa+1} |\xi|^{1/2} \\
    &\bigg(\cos(|\xi|-\frac{\rho\pi}2-\frac\pi4)\J_{1-\rho}(|\xi|)-\cos(|\xi|+\frac{\rho\pi}2-\frac\pi4)\J_{\rho-1}(|\xi|)\bigg)\\
     &+\mathcal O((1+t)^{-1}|\xi|^{-1})\\
\end{align*}
as $t\to\infty$, $|\xi|\ge c$. Hence, as $t\to\infty$ the function $m_{1,1}(t,\xi)$ tends uniformly (on $|\xi|\ge c$) to
the limit $m_{1,1}(\infty,\xi)$ which extends continuously (and analytically) up to $\xi=0$. The last statement follows from
the generalized power series expansion of $J_\rho(z)=z^{\rho}\Lambda_\rho(z)$, $\Lambda_\rho(0)\ne0$ together with
$\kappa+3/2=(\mu+1)/2=1-\rho$.

A similar calculation yields for the other entries of the matrix limit expressions of the form
\begin{align*}
  m_{1,2}(\infty,\xi)&=\sqrt{\frac\pi2}\csc(\rho\pi) [\xi]^{\kappa} |\xi|^{1/2} \\
    &\bigg(\cos(|\xi|-\frac{\rho\pi}2-\frac\pi4)\J_{-\rho}(|\xi|)+\cos(|\xi|+\frac{\rho\pi}2-\frac\pi4)\J_{\rho}(|\xi|)\bigg)\\
  m_{2,1}(\infty,\xi)&=\sqrt{\frac\pi2}\csc(\rho\pi) [\xi]^{\kappa+1} |\xi|^{1/2} \\
    &\bigg(\sin(|\xi|-\frac{\rho\pi}2-\frac\pi4)\J_{1-\rho}(|\xi|)+\sin(|\xi|+\frac{\rho\pi}2-\frac\pi4)\J_{\rho-1}(|\xi|)\bigg)\\
  m_{2,2}(\infty,\xi)&=\sqrt{\frac\pi2}\csc(\rho\pi) [\xi]^{\kappa} |\xi|^{1/2} \\
    &\bigg(\sin(|\xi|-\frac{\rho\pi}2-\frac\pi4)\J_{-\rho}(|\xi|)-\sin(|\xi|+\frac{\rho\pi}2-\frac\pi4)\J_{\rho}(|\xi|)\bigg).
\end{align*}
By analytic continuation these formulas are valid for all $\xi\in\R^n$ and $m_{i,j}(\infty,\cdot)\in L^\infty(\R^n)$. Furthermore,
all entries are non-zero for $\xi=0$ (due to the exact cancellation).

Furthermore the representations are analytic in $\rho$ using the definition of the Weber functions 
$\mathcal Y_\rho(z)=\csc(\rho\pi)\big(\J_\rho(z)\cos(\rho\pi)-\J_{-\rho}(z))$. 
\end{proof}

An application of Liouville theorem to the differential equation for $(|\D|u,u_t)^T$ yields
\begin{equation}
 \begin{vmatrix} \Psi_{1,0,\rho-1,1} & -\Psi_{1,0,\rho,0}\\
  \Psi_{1,0,\rho-1,0} & -\Psi_{1,0,\rho,-1} \end{vmatrix} = (1+t)^{-\mu}
\end{equation}
and therefore $\det (1+t)^{\mu/2} \mathbb E^{(\kappa)}(t,\xi)=[\xi]^{1+2\kappa}$.

\begin{cor}
  The operator $Z_+:[\D]^\kappa L^2(\R^n)\to L^2(\R^n)$ is injective and its symbol satisfies $\det Z_+(\xi)=[\xi]^{1+2\kappa}$.
\end{cor}

Remark. What have we obtained so far? The existence of nontrivial Cauchy data $(\langle\D\rangle u_1,u_2)$ in the null space of 
$Z_+$ is equivalent to the fact that the energy of the corresponding solution 
$||\mathbb E(t)(\langle\D\rangle u_1,u_2)^T||_2^2$ decays faster than $t^{-\mu}$. But we have shown, that such Cauchy data do not
exist. This means, we have proven a two-sided energy estimate.

\begin{cor}
  Let $(\langle\D\rangle u_1,u_2)\in[\D]^\kappa L^2(\R^n)$. Then 
  \begin{equation}
    E(u;t)\sim (1+t)^{-\mu}.
  \end{equation} 
\end{cor} 

Furthermore, the solutions behave in energy space as solutions to the free wave equation multiplied by the decay rate 
$(1+t)^{-\mu/2}$. In the case of high frequencies we can even say more. The operator $Z_+$ almost preserves high frequencies. Using again the asymptotic representation of Bessel functions for large arguments, Proposition~\ref{pB},
we conclude

\begin{cor}
  It holds $\lim_{|\xi|\to\infty} Z_+(\xi)=I$.
\end{cor}

For $\mu=2$ we have even more. In this case the Bessel functions are trigonometric ones and $Z_+(\xi)=I$.

In order to conclude this article we will give one further interpretation to the assumptions on the Cauchy data
we made.
In case $0\le\kappa<\frac n2$ the Sobolev-Hardy inequality can be used to obtain an embedding of the space of
weighted $L^2$-functions into $[\D]^\kappa L^2$.

\begin{lem}
  Let $0\le\kappa<\frac n2$. Then $\langle x\rangle^{-\kappa}L^2(\R^n)\subseteq [\D]^\kappa L^2(\R^n)$.
\end{lem}
\begin{proof}
Let $\chi\in C_0^\infty(\R^n)$ satisfy $\chi(\xi)=1$ near $\xi=0$. Then
\begin{multline*}
  ||f||_{[\D]^\kappa L^2}= ||[\xi]^{-\kappa} \hat f||_2\sim ||\,|\xi|^{-\kappa}\chi(\xi)\hat f||_2+||(1-\chi(\xi))f||_2\\
    \lesssim ||\hat f||_{H^\kappa} + ||\hat f||_2\sim ||\langle x\rangle^\kappa f||_2 
\end{multline*}
by Sobolev-Hardy inequality and Plancherel's theorem.
\end{proof}

\end{document}